\documentclass{amsart}
\usepackage{amsmath, amssymb}
\usepackage{amsfonts}
\usepackage{amsthm}
\usepackage[arrow,matrix,curve,cmtip,ps]{xy}

\usepackage{amsthm}

\allowdisplaybreaks

\newtheorem{theorem}{Theorem}[section]
\newtheorem{lemma}[theorem]{Lemma}
\newtheorem{proposition}[theorem]{Proposition}

\newtheorem*{theorem*}{Theorem}
\theoremstyle{remark}

\newtheorem{definition}[theorem]{Definition}

\numberwithin{equation}{section}

\newcommand{\Z}{\mathbb{Z}}
\newcommand{\N}{\mathbb{N}}
\newcommand{\R}{\mathbb{R}}
\newcommand{\C}{\mathbb{C}}

\newcommand{\K}{\mathcal{K}}

\newcommand{\p}{\mathcal{P}}

\newcommand{\coker}{\operatorname{coker}}

\renewcommand{\i}{\iota}

\begin{document}
\title[The ordered $K_0$-group of a graph $C^*$-algebra]{The
ordered $\boldsymbol{K_0}$-group of a graph $\boldsymbol{C^*}$-algebra}
\author{Mark Tomforde 
}

\address{Department of Mathematics\\ Dartmouth College\\
Hanover\\ NH 03755\\ USA}

\curraddr{Department of Mathematics\\ University of Iowa\\
Iowa City\\ IA 52242\\ USA}

\email{tomforde@math.uiowa.edu}


\date{\today}
\subjclass{46L55}

\begin{abstract}
We calculate the ordered $K_0$-group of a graph $C^*$-algebra and mention
applications of this result to AF-algebras, states on the $K_0$-group of a
graph algebra, and tracial states of graph algebras.
\end{abstract}

\maketitle

\section{Preliminaries} \label{Pre}

We provide some basic facts about graph algebras and refer
the reader to \cite{KPR}, \cite{BPRS}, and \cite{RS} for more
details.  A (directed) graph $E=(E^0, E^1, r, s)$ consists of a
countable set $E^0$ of vertices, a countable set $E^1$ of edges,
and maps $r,s:E^1 \rightarrow E^0$ identifying the range and
source of each edge.  A vertex $v \in E^0$ is called a
\emph{sink} if $|s^{-1}(v)|=0$, and $v$ is called an
\emph{infinite emitter} if $|s^{-1}(v)|=\infty$.  If $v$ is either a sink or an
infinite emitter, we call $v$ a \emph{singular vertex}.  A graph
$E$ is said to be \emph{row-finite} if it has no
infinite emitters.  The \emph{vertex matrix} of $E$ is the square matrix $A$
indexed by the vertices of $E$ with $A(v,w)$ equal to the number of edges from
$v$ to $w$.

If $E$ is a graph we define a \emph{Cuntz-Krieger $E$-family}
to be a set of mutually orthogonal projections $\{p_v : v \in
E^0\}$ and a set of partial isometries $\{s_e : e \in E^1\}$ with
orthogonal ranges which satisfy the \emph{Cuntz-Krieger relations}:
\begin{enumerate}
\item $s_e^* s_e = p_{r(e)}$ for every $e \in E^1$;
\item $s_e s_e^* \leq p_{s(e)}$ for every $e \in E^1$;
\item $p_v = \sum_{\{e : s(e)=v\}} s_e s_e^*$ for every $v
\in E^0$ that is not a singular vertex.
\end{enumerate}
The \emph{graph algebra $C^*(E)$} is defined to be the
$C^*$-algebra generated by a universal Cuntz-Krieger $E$-family.

The graph algebra $C^*(E)$ is unital if and only if $E$ has a
finite number of vertices, cf. \cite[Proposition~1.4]{KPR}, and in
this case $1_{C^*(E)} = \sum_{v \in E^0} p_v$.  If $E$ has an
infinite number of vertices, and we list the vertices of $E$ as
$E^0 = \{v_1, v_2, \ldots \}$ and define $p_n := \sum_{i=1}^n
p_{v_i}$, then $\{ p_n \}_{n=1}^\infty$ will be an approximate
unit for $C^*(E)$.

\section{The ordered $K_0$-group}

If $A$ is a $C^*$-algebra let $\p (A)$ denote the set
of projections in $A$.  It is a fact that if $A$ is unital
(or more generally, if $A \otimes \K$ admits an approximate unit
consisting of projections), then $K_0(A) = \{ [p]_0-[q]_0 : p,q
\in \p (A \otimes\K) \}$.  In addition, the positive cone $K_0(A)^+
= \{ [p]_0 : p \in \p (A \otimes \K) \}$ makes $K_0(A)$ a
pre-ordered abelian group.  If $A$ is also stably finite, then
$(K_0(A),K_0(A)^+)$ will be an ordered abelian group.

Here we compute the positive cone of the $K_0$-group of a graph
$C^*$-algebra.  Throughout this section we let $\Z^K$ and $\N^K$ denote
$\bigoplus_K \Z$ and $\bigoplus_K \N$, respectively.

\begin{lemma} \label{K-theory-row-finite}
Let $E = (E^0,E^1,r,s)$ be a row-finite graph.  Also let $W$ denote the set
of sinks of $E$ and let $V:=E^0 \backslash W$.  Then
with respect to the decomposition $E^0 = V \cup W$ the vertex
matrix of $E$ will have the form $$A_E = \begin{pmatrix} B & C \\
0 & 0 \end{pmatrix}.$$  For $v \in E^0$, let $\delta_v$ denote the
element of $\Z^V \oplus \Z^W$ with a 1 in the $v^{\text{th}}$ entry and
$0$'s elsewhere. 

If we consider $\begin{pmatrix} B^t-I \\ C^t \end{pmatrix} :
\Z^V \rightarrow \Z^V \oplus \Z^W$, then  $K_0(C^*(E)) \cong
\coker \begin{pmatrix} B^t-I \\ C^t \end{pmatrix}$ via an
isomorphism which takes $[p_v]_0$ to $[\delta_v]$ for each $v
\in E^0$.  Furthermore, this isomorphism takes $(K_0(C^*(E)))^+$
to $\{ [x] : x \in \N^V \oplus \N^W \}$, where $[x]$ denotes the
class of $x$ in $\coker \begin{pmatrix} B^t-I \\ C^t
\end{pmatrix}$.
\end{lemma}

\begin{proof}
The fact that $K_0(C^*(E)) \cong
\coker \begin{pmatrix} B^t-I \\ C^t \end{pmatrix}$ is shown for
row-finite graphs in \cite[Theorem~3.1]{RS}.  Thus all that remains to be
verified in our claim is that this isomorphism identifies
$(K_0(C^*(E)))^+$ with $\{ [x] : x \in \N^V \oplus \N^W \}$.  To do
this, we will simply examine the proof of \cite[Theorem~3.1]{RS} to
determine how the isomorphism acts.  We will assume that the reader is
familiar with this proof, and use the notation established in it
without comment.

If $E \times_1 [m,n]$ is the graph defined in \cite[Theorem~3.1]{RS}, then we
see that $E \times_1 [m,n]$ is a row-finite graph with no loops and in which
every path has length at most $n-m$.  Therefore we can use the arguments in the
proofs of \cite[Proposition~2.1]{KPR}, \cite[Corollary~2.2]{KPR}, and
\cite[Corollary~2.3]{KPR} to conclude that $C^*(E \times_1 [m,n])$ is the
direct sum of copies of the compact operators (on spaces of varying
dimensions), indexed by the sinks of $E \times_1 [m,n]$ and that each summand
contains precisely one projection $p_{(v,k)}$ associated to a sink as a minimal
projection.  Thus $$K_0(C^*(E \times_1 [m,n])) \cong \bigoplus_{v \in V}
\Z [p_{(v,n)}]_0 \oplus \bigoplus_{k=0}^{n-m} \bigoplus_{v \in W} \Z
[p_{(v,n-k)}]_0$$ and $K_0(C^*(E \times_1
[m,n]))^+$ is identified with  $$\bigoplus_{v \in V} \N [p_{(v,n)}]_0 \oplus
\bigoplus_{k=0}^{n-m} \bigoplus_{v \in W} \N [p_{(v,n-k)}]_0.$$
By the continuity of $K$-theory, one can let $m$ tend to $-\infty$ and
deduce that 
\begin{align*}
K_0(C^*(E \times_1 [-\infty,n])) &\cong \bigoplus_{v \in V} \Z
[p_{(v,n)}]_0 \oplus \bigoplus_{k=0}^{\infty} \bigoplus_{v \in W}
\Z [p_{(v,n-k)}]_0 \\
&\cong \Z^V \oplus Z^W \oplus Z^W \oplus \ldots.
\end{align*}
Furthermore, it follows from \cite[Theorem~6.3.2(ii)]{RLL} that
this isomorphism identifies $K_0(C^*(E \times_1 [-\infty,n]))^+$
with $\N^V \oplus \N^W \oplus \N^W \oplus \ldots$.

This computation is used later in the proof of
\cite[Theorem~3.1]{RS}, where the $K_0$ functor is applied to a
commutative diagram to obtain Figure~(3.5) of \cite{RS}, which we
reproduce here:
\begin{equation*}
\xymatrix{ \Z^V \oplus \Z^W \oplus \Z^W \oplus \ldots \ar[d]_{1-D}
\ar[r]^D & \Z^V \oplus \Z^W \oplus \Z^W \oplus \ldots
\ar[d]_{1-D} \ar[r]^-{\iota_*^{n+1}} & K_0(C^*(E \times_1 \Z))
\ar[d]_{1-\beta_*^{-1}} \\ 
\Z^V \oplus \Z^W \oplus \Z^W \oplus \ldots \ar[r]^{D} & \Z^V
\oplus \Z^W \oplus \Z^W \oplus \ldots \ar[r]^-{\iota_*^{n+1}} &
K_0(C^*(E \times_1 \Z)) }
\end{equation*}

Now it is shown in \cite[Lemma~3.3]{RS} that the homomorphism
$\iota_*^1$ induces an isomorphism $\overline{\iota}_*^1$ of
$\coker (1-D)$ onto $\coker (1-\beta_*^{-1}) = K_0(C^*(E))$.  We
shall show that $\overline{\iota}_*^1 ( \N^V \oplus \N^W \oplus
\N^W \oplus \ldots ) = K_0(C^*(E))^+$.  To begin, note that it
follows from \cite[Theorem~6.3.2(ii)]{RLL} that $$K_0(C^*(E
\times_1 \Z))^+ = \bigcup_{n=1}^\infty \iota_*^n ( \N^V \oplus
\N^W \oplus \N^W \oplus \ldots ).$$  Since $\coker
(1-\beta_*^{-1}) = K_0(C^*(E))$, this implies that $$K_0(C^*(E))^+
= \bigcup_{n=1}^\infty \{ [\iota_*^n(y)] : y \in \N^V \oplus \N^W
\oplus \N^W \oplus \ldots \}$$ where $[\iota_*^n(y)]$ denotes the
equivalence class of $\iota_*^n(y)$ in $\coker (1-\beta_*^{-1})$. 
We shall show that the right hand side of this equation is equal
to $\{ [\iota_*^1(y)] : y \in \N^V \oplus \N^W \oplus \N^W \oplus
\ldots \}$.  Let $[\iota_*^n(y)]$ be a typical element in the right
hand side.  Then from the commutativity of the above diagram
$\iota_*^n(y) - \iota_*^n(Dy) = \iota_*^n((1-D)y) =
(1-\beta_*^{-1}) (\iota_*^n(y))$ which is $0$ in $\coker
(1-\beta_*^{-1})$.  But then $\iota_*^1(y) = \iota_*^n(D^{n-1}y) =
\iota_*^n(y)$ in $\coker (1-\beta_*^{-1})$.  Hence 
\begin{equation} \label{iota-map}
K_0(C^*(E))^+ = \{ [\iota_*^1(y)] : y \in \N^V \oplus \N^W \oplus
\N^W \oplus \ldots \}. 
\end{equation}

Next, recall that \cite[Lemma~3.4]{RS} shows that the inclusion $j
: \Z^V \oplus Z^W \hookrightarrow \Z^V \oplus \Z^W \oplus \Z^W
\oplus \ldots$ induces an isomorphism $\overline{j}$ of $\coker K$
onto $\coker (1-D)$. We wish to show that 
\begin{equation} \label{j-map}
\overline{j} (\{[x] : x \in \N^V \oplus \N^W \} ) = \{ [y] : y \in
\N^V \oplus \N^W \oplus \N^W \oplus \ldots \}.
\end{equation}
It suffices to show that any element $(n,m_1,m_2,m_3, \ldots) \in
\N^V \oplus \N^W \oplus \N^W \oplus \ldots$ is equal to an element
of the form $(a, b, 0, 0, 0, \ldots)$ in $\coker(1-D)$.  But given
$(n,m_1,m_2,m_3, \ldots)$ we see that since this element is in the
direct sum, there exists a positive integer $k$ for which $i>k$
implies $m_i=0$.  Thus
$$ \begin{pmatrix} n \\ m_1 + \ldots +m_k \\ 0 \\
\vdots \\ 0 \\ 0 \\ \vdots\end{pmatrix} - \begin{pmatrix} n \\ m_1 \\ m_2 \\ \vdots \\ m_k \\
0 \\ \vdots\end{pmatrix} =  \begin{pmatrix} 1-B^t & 0
& 0 & 0 & .\\ -C^t & 1 & 0 & 0 & . \\ 0 & -1 & 1 & 0 & . \\ 0 &
0 & -1 & 1 & . \\ . & . & . & . & \ddots \end{pmatrix} 
\begin{pmatrix} 0 \\ m_2 + \ldots +m_k \\ m_3 + \ldots + m_k \\
m_4 + \ldots + m_k \\ \vdots \\ m_k \\ 0 \\ \vdots\end{pmatrix}
$$ and so $(n,m_1,m_2, \ldots)$ equals $(n,m_1+\ldots + m_k, 0, 0,
\ldots)$ in $\coker (1-D)$, and (\ref{j-map}) holds.

Finally, the isomorphism between $\coker \begin{pmatrix}
B^t-I \\ C^t \end{pmatrix}$ and $K_0(C^*(E))$ is defined to be
$\overline{\iota}_*^1 \circ \overline{j}$.  But (\ref{j-map}) and
(\ref{iota-map}) show that this isomorphism takes $\{ [x] : x \in
\N^V \oplus \N^W \}$ onto $(K_0(C^*(E)))^+$.
\end{proof}

\begin{theorem} \label{K-theory}
Let $E = (E^0,E^1,r,s)$ be a graph.  Also let $W$ denote the set
of singular vertices of $E$ and let $V:=E^0 \backslash W$.  Then
with respect to the decomposition $E^0 = V \cup W$ the vertex
matrix of $E$ will have the form $$A_E = \begin{pmatrix} B & C \\
* & * \end{pmatrix}$$ where $B$ and $C$ have entries in
$\Z$ and the $*$'s have entries in $\Z \cup \{ \infty \}$.  Also
for $v \in E^0$, let $\delta_v$ denote the element of
$\Z^V \oplus \Z^W$ with a 1 in the $v^{\text{th}}$ entry and $0$'s
elsewhere. 

If we consider $\begin{pmatrix} B^t-I \\ C^t \end{pmatrix} :
\Z^V \rightarrow \Z^V \oplus \Z^W$, then  $K_0(C^*(E)) \cong
\coker \begin{pmatrix} B^t-I \\ C^t \end{pmatrix}$ via an
isomorphism which takes $[p_v]_0$ to $[\delta_v]$ for each $v
\in E^0$.  Furthermore, this isomorphism takes $(K_0(C^*(E)))^+$
onto the semigroup generated by $\{ [\delta_v] : v \in E^0 \} \cup \{
[\delta_v] - \sum_{e \in S}[\delta_{r(e)}] : \text{$v$ is an infinite
emitter and $S$ is a finite subset of $s^{-1}(v)$} \}$.
\end{theorem}

\begin{proof}
The fact that $K_0(C^*(E)) \cong \coker \begin{pmatrix} B^t-I \\ C^t
\end{pmatrix}$ was established in \cite[Theorem~3.1]{DT2} using the
isomorphisms constructed in \cite[Lemma~2.3]{DT2}.  We shall examine the proof
of \cite[Theorem~3.1]{DT2} to determine where the positive cone of
$K_0(C^*(E))$ is sent.  Again, we shall assume that the reader is
familiar with the proof, and use the notation established in it without
comment.

We begin by letting $F$ denote a desingularization of $E$ (see
\cite[\S2]{DT2}).  Then \cite[Theorem~2.11]{DT1} shows that there exists a
homomorphism $\phi :C^*(E) \to C^*(F)$ which embeds $C^*(E)$ onto a full corner
of $C^*(F)$ and takes each $p_v$ to the projection in $C^*(F)$ corresponding
to $v$.  Since $\phi$ is an embedding onto a full corner, it induces an
isomorphism $\phi_* : K_0(C^*(E)) \to K_0(C^*(F))$ which takes the class of
$p_v$ in $K_0(C^*(E))$ to the class of the corresponding projection in
$K_0(C^*(F))$.  By Theorem~\ref{K-theory-row-finite}, if
$A_F$ denotes the vertex matrix of $F$, then $K_0(C^*(E)) \cong \coker
(A_F^t-I)$ and $K_0(C^*(E))^+$ is identified with $\{ [x] : x \in \bigoplus_V
\Z \oplus \bigoplus_W \Z \oplus \bigoplus_W Q \}$ where $Q := \bigoplus_\N
\Z$.  Now it is shown in the proof of \cite[Lemma~2.3]{DT2} that the inclusion
map $\rho : \bigoplus_V \Z \oplus \bigoplus_W \Z \to \bigoplus_V \Z \oplus
\bigoplus_W \Z \oplus \bigoplus_W Q$ induces an isomorphism $\overline{\rho} : 
\coker \begin{pmatrix} B^t-I \\ C^t \end{pmatrix} \to \coker(A_F-I)$.  Since
this isomorphism identifies the class of $\delta_v \in \bigoplus_V \Z \oplus
\bigoplus_W \Z$ with the class of $\begin{pmatrix} \delta_v \\ 0
\end{pmatrix} \in \bigoplus_V \Z \oplus \bigoplus_W \Z \oplus 
\bigoplus_W Q$, it follows that $[p_v]_0 \in K_0(C^*(E))$ is identified
with $[\delta_v] \in \coker \begin{pmatrix} B^t-I \\ C^t
\end{pmatrix}$.  

All that remains is to determine where this isomorphism sends the positive
cone of $K_0(C^*(E))$.  Let $\Gamma$ denote the semigroup of elements that
$\overline{\rho}$ sends to $\{ [x] : x \in \bigoplus_V \Z \oplus \bigoplus_W \Z
\oplus \bigoplus_W Q \}$.  Now certainly $\{ [\delta_v] : v \in E^0 \}$ is in
$\Gamma$.  Furthermore, for any infinite emitter $v$ and finite subset $S
\subseteq s^{-1}(v)$ we have that $$[p_v]_0 - \sum_{e \in S} [p_{r(e)}]_0 =
[p_v]_0 - \sum_{e \in S} [s_e^*s_e]_0 = [p_v]_0 - \sum_{e \in S} [s_es_e^*]_0
= [p_v - \sum_{e \in S} s_es_e^*]_0$$ and this element belongs to
$K_0(C^*(E))^+$.  Since $K_0(C^*(E))^+$ is identified with $\{ [x] : x \in
\bigoplus_V \Z \oplus \bigoplus_W \Z \oplus \bigoplus_W Q \}$ this implies that
the class of $ \begin{pmatrix} \delta_v \\ 0 \end{pmatrix}  - \sum_{e \in S}
\begin{pmatrix} \delta_{r(e)} \\ 0 \end{pmatrix}$ is in $\{ [x] : x \in
\bigoplus_V \Z \oplus \bigoplus_W \Z \oplus \bigoplus_W Q \}$ and thus $[p_v] -
\sum_{e \in S} [p_{r(e)}]$ is in $\Gamma$.  On the other hand, we know that
$\Gamma$ is generated by the elements that $\overline{\rho}$ sends to the
classes of the generators of $\bigoplus_V \Z \oplus \bigoplus_W \Z \oplus
\bigoplus_W Q$.  Now certainly the inverse image under $\overline{\rho}$ of
the class of $\begin{pmatrix} \delta_v \\ 0 \end{pmatrix}$ for $v \in V \cup W$
will be $[\delta_v]$.  In addition, if $v_i$ is a vertex on the tail added to
an infinite emitter $v$, then we see that the inverse image under
$\overline{\rho}$ of the element $\begin{pmatrix} 0 \\ \delta_{v_i} 
\end{pmatrix}$ will be $\begin{pmatrix} \mathbf{u} \\ \mathbf{v}
\end{pmatrix}$ where $\mathbf{u}$ and $\mathbf{v}$ are as defined in the final
paragraph of \cite[Lemma~2.3]{DT2}.  However, one can verify from how
$\mathbf{u}$ and $\mathbf{v}$ are defined that $\begin{pmatrix} \mathbf{u} \\
\mathbf{v} \end{pmatrix}$ will have the form $\delta_v - \sum_{e \in S}
\delta_{r(e)}$ for some finite $S \subseteq s^{-1}(v)$.  Thus $\Gamma$ is
generated by the elements $[\delta_v]$ and $[\delta_v] - \sum_{e \in S}
[\delta_{r(e)}]$.
\end{proof}

\section{Applications}

\subsection{AF-algebras}

The graph algebra $C^*(E)$ is an AF-algebra if and only if $E$ has no loops
\cite[Theorem~2.4]{KPR}.  By Elliott's Theorem
AF-algebras are classified by their ordered $K_0$-groups.  Hence for two graphs
containing no loops, Theorem~\ref{K-theory} can be used to determine if
their associated $C^*$-algebras are isomorphic (as well as stably isomorphic).

\subsection{States on $\boldsymbol{K_0(C^*(E))}$}

If $A$ is a $C^*$-algebra containing a countable approximate unit $\{ p_n
\}_{n=1}^\infty$ consisting of projections, then a \emph{state} on $K_0(A)$ is
a homomorphism $f : K_0(A) \to \R$ such that $f(K_0(A)^+) \subseteq \R^+$ and
$\lim_{n \to \infty} f([p_n]_0) = 1$.  The set of all states on $K_0(A)$ is
denoted $S(K_0(A))$ and we make it into a topological space by giving it the
weak-$*$ topology.

\begin{definition} \label{graph-trace}
If $E$ is a graph, then a \emph{graph trace} on $E$ is a
function $g : E^0 \rightarrow \R^+$ with the
following two properties: 
\begin{enumerate}
\item \label{g-t-1} For any nonsingular vertex $v \in E^0$ we have $g(v) =
\sum_{ \{e \in E^1 : s(e) = v \} } g(r(e)).$
\item \label{g-t-2} For any infinite emitter $v \in E^0$ and
any finite set of edges $e_1, \ldots, e_n \in s^{-1}(v)$ we have $g(v) \geq
\sum_{i=1}^n g(r(e_i)).$
\end{enumerate}
We define the norm of $g$ to be the (possibly infinite) value $\|g
\| := \sum_{v \in E^0} g(v)$, and we shall use $T(E)$ to denote the set of all
graph traces on $E$ with norm 1.
\end{definition}

\begin{proposition} \label{g-t-is-state}
If $E$ is a graph, then the state space $S(K_0(C^*(E)))$ with the weak-$*$
topology is naturally isomorphic to $T(E)$ with the topology generated
by the subbasis $\{ N_{v,\epsilon} (g) : v \in E^0, \epsilon > 0,
\text{ and } g \in T(E) \}$, where $N_{v,\epsilon} (g):= \{h \in
T(E) : |h(v)-g(v)| < \epsilon \}$.
\end{proposition}

\begin{proof}
We define a map $\i : S(K_0(C^*(E))) \to T(E)$ by
$\i(f)(v) := f([p_v]_0)$.  We shall show that $\i$ is an affine
homeomorphism.  To see that $\i$ is injective note that if $\i(f_1)=\i(f_2)$, 
then for each $v \in E^0$ we have that $f_1([p_v]_0) = \i (f_1) (v) =
\i (f_2) (v) = f_2([p_v]_0)$, and since the $[p_v]_0$'s generate
$K_0(C^*(E))$ it follows that $f_1 = f_2$.

To see that $\i$ is surjective, let $g : E^0 \to \R^+$
be a graph trace.  We shall define a homomorphism $f : \coker \begin{pmatrix}
B^t-I \\ C^t \end{pmatrix} \to \R$ by setting $f ([\delta_v]) := g(v)$. 
Because $g$ satisfies (\ref{g-t-1}) of Definition~\ref{graph-trace} we see
that $f$ is well defined.  Also, since the values of $g$ are positive and $g$
satisfies (\ref{g-t-2}) of Definition~\ref{graph-trace} we see that
$f(K_0(C^*(E))^+) \subseteq \R^+$.  Finally, since $g$ has norm 1 we see that
$\lim_{n \to \infty} f ( \sum_{i=1}^n p_{v_i} ) = \lim_{n \to \infty}
\sum_{i=1}^n g (v_i) = \| g \| =1$.  So $f$ is a state on $K_0(C^*(E))$ and $\i
(f) = g$.

It is straightforward to verify that $\i$ is an affine homeomorphism.
\end{proof}

\subsection{Tracial states on $\boldsymbol{C^*(E)}$}

A \emph{trace} on a $C^*$-algebra $A$ is a linear functional $\tau : A
\rightarrow \C$ with the property that $\tau(ab)=\tau(ba)$ for all $a,b \in
A$.  We say that $\tau$ is \emph{positive} if $\tau (a) \geq 0$ for all $a
\in A^+$.  If $\tau$ is positive and $\| \tau \| = 1$ we call
$\tau$ a \emph{tracial state}.  The set of all tracial states is
denoted $T(A)$ and when $T(A)$ is nonempty we equip it with the
weak-$*$ topology.  Let $A$ be a $C^*$-algebra with a countable approximate
unit consisting of projections.  If $\tau$ is a trace on $A$,
then it induces a map $K_0(\tau) : K_0(A) \to \R$ given by
$K_0(\tau)([p]_0-[q]_0) = \tau(p)-\tau(q)$.  The map $K_0(\tau)$ will be an
element of $S(K_0(A))$ (see \cite[\S 5.2]{RLL} for more details) and thus there
is a continuous affine map $r_A : T(A) \to S(K_0(A))$ defined by $r_A(\tau) :=
K_0(\tau)$.  

It is a fact that any quasitrace on an exact $C^*$-algebra extends to a trace
(this was proven by Haagerup for unital $C^*$-algebras \cite{Haa} and shown to
hold for nonunital $C^*$-algebras by Kirchberg \cite{Kir2}).  Furthermore,
Blackadar and R\o rdam showed in \cite{BR} that when $A$ is unital every
element in $K_0(A)$ lifts to a quasitrace.  It is straightforward to extend the
result of Blackadar and R\o rdam to $C^*$-algebras with a countable
approximate unit consisting of projections.  Thus when $A$ is a graph algebra
we see that the map $r_A : T(A) \to S(K_0(A))$ is surjective.

If $A$ has real rank zero, then the span of
the projections in $A$ is dense in $A$ and  $r_A$ is
injective.  It was shown in \cite{Jeo} that a graph
algebra $C^*(E)$ has real rank zero if and only if the graph $E$ satisfies
Condition~(K); that is, no vertex in $E$ is the base of exactly one simple
loop.  Therefore, when $A = C^*(E)$ and $E$ is a graph satisfying
Condition~(K), the map $r_A$ is a homeomorphism and
Proposition~\ref{g-t-is-state} shows that the tracial states on $C^*(E)$ are
identified in a canonical way with $T(E)$.


\begin{thebibliography}{99}

\bibitem{BPRS}
T.~Bates, D.~Pask, I.~Raeburn, and W.~Szyma\'nski,
\emph{The $C^*$-algebras of row-finite graphs}, New York J.
Math. \textbf{6} (2000), 307--324.

\bibitem{BR}
B.~Blackadar and M.~R\o rdam, \emph{Extending states on 
preordered semigroups and the existence of quasitraces on 
$C^*$-algebras}, J.~Algebra, \textbf{152} (1992), 240--247.

\bibitem{DT1}
D.~Drinen and M.~Tomforde, \emph{The $C^*$-algebras of
arbitrary graphs}, Rocky Mountain J. Math, to appear.

\bibitem{DT2}
D.~Drinen and M.~Tomforde, \emph{Computing $K$-theory and
Ext for graph $C^*$-algebras}, Illinois J. Math., to
appear.

\bibitem{Haa}
U.~Haagerup, \emph{Every quasi-trace on an exact $C^*$-algebra is 
a trace}, preprint (1991).

\bibitem{Jeo}
J.~A~Jeong, \emph{Real rank of generalized Cuntz-Krieger algebras}, 
preprint (2000).

\bibitem{Kir2}
E.~Kirchberg, \emph{On the existence of traces on exact stably
projectionless simple $C^*$-algebras},  Operator Algebras and
their Applications (P.~A.~Fillmore and J.~A.~Mingo, eds.) Fields
Institute Communications, vol. 13, Amer. Math. Soc. 1995, p.
171--172.

\bibitem{KPR}
A.~Kumjian, D.~Pask, and I.~Raeburn, \emph{Cuntz-Krieger
algebras of directed graphs}, Pacific J. Math. \textbf{184}
(1998), 161--174.

\bibitem{RS}
I.~Raeburn and W.~Szyma\'nski,
\emph{Cuntz-Krieger algebras of infinite graphs and matrices},
preprint (2000).

\bibitem{RLL}
M.~R\o rdam, F.~Larsen, and N.~J.~Laustsen, An Introduction
to $K$-theory for $C^*$-algebras, London Mathematical
Society Student Texts \textbf{49}, Cambridge University
Press, Cambridge, 2000.

\end{thebibliography}
\end{document}